\documentclass[11pt]{article}
\usepackage{amsfonts,amsmath,amssymb,tikz}
\DeclareMathSymbol{\lcurs}{\mathord}{letters}{"0B}
\newtheorem{theorem}{\bf Theorem}[section]
\newtheorem{corollary}[theorem]{\bf Corollary}
\newtheorem{lemma}[theorem]{\bf Lemma}
\newtheorem{proposition}[theorem]{\bf Proposition}

\newtheorem{problem}[theorem]{\bf Problem}

\newcommand{\proof}{\noindent{\bf Proof.\ }}
\newcommand{\qed}{\hfill $\square$ \bigskip}
\newcommand{\cp}{\,\square\,}

\begin{document}

\title{Asymptotic properties of Fibonacci cubes and Lucas cubes}

\author{
Sandi Klav\v zar
\\
Faculty of Mathematics and Physics \\
University of Ljubljana, Slovenia \\
and \\
Faculty of Natural Sciences and Mathematics \\
University of Maribor, Slovenia \\
and \\
Institute of Mathematics, Physics and Mechanics, Ljubljana \\
\tt{sandi.klavzar@fmf.uni-lj.si} 
\and
Michel Mollard \\
CNRS Universit\'e Joseph Fourier\\
Institut Fourier, BP 74 \\
100 rue des Maths, 38402 St Martin d'H\`eres Cedex, France \\
\tt{michel.mollard@ujf-grenoble.fr}
}

\date{}

\maketitle

\begin{abstract}
It is proved that the asymptotic average eccentricity and the asymptotic average degree of both Fibonacci cubes and Lucas cubes are $(5+\sqrt 5)/10$ and $(5-\sqrt 5)/5$, respectively. A new labeling of the leaves of Fibonacci trees is introduced and it is proved that the eccentricity of a vertex of a given Fibonacci cube is equal to the depth of the associated leaf in the corresponding Fibonacci tree. Hypercube density is also introduced and studied. The hypercube density of both Fibonacci cubes and Lucas cubes is shown to be $(1-1/\sqrt 5)/\log_2\varphi$, where $\varphi$ is the golden ratio, and the Cartesian product of graphs is used to construct families of graphs with a fixed, non-zero hypercube density. It is also proved that the average ratio of the numbers of Fibonacci strings with a 0 resp. a 1 in a given position, where the average is taken over all positions, converges to $\varphi^2$, and likewise for Lucas strings.
\end{abstract}


\medskip\noindent
{\bf AMS Subj. Class.:} 05C12, 40A05, 05A15

\section{Introduction}
\label{sec:intro}

Fibonacci cubes~\cite{hsu-1993} and Lucas cubes~\cite{mupe-2001} form appealing infinite families of graphs which are the focus of much current research; see the recent survey~\cite{kl-2013}. These cubes are subgraphs of hypercubes and on one hand, they inherit many of the fine properties of hypercubes, while on the other hand their size grows significantly slower than that of hypercubes. Moreover, Fibonacci cubes and Lucas cubes found several applications, for instance in theoretical chemistry; see~\cite{zhou-2009} for a use of Fibonacci cubes and~\cite{zibe-2012} for a use of Lucas cubes. We also mention that in~\cite{ra-2013} an investigation of the on-line routing of linear permutations on these cubes is performed.

In the last years large graphs (and/or complex networks) became a topic of great interest---not only in mathematics but also elsewhere---hence it seems justified to consider the asymptotic behaviour of applicable families of graphs, such as Fibonacci and Lucas cubes. The average normed distance of these graphs was  proved in~\cite{klmo-2012} to be $2/5$. In this paper we study the limit behaviour of the average eccentricity of these cubes, the limit behaviour of their average degree, and some related topics. For some general properties of the average eccentricity see~\cite{dago-2004,il-2012}, while~\cite{hipa-2012} gives the the average eccentricity of Sierpi\'nski graphs. 

We proceed as follows. In the rest of this section concepts needed in this paper are formally introduced. In Section~\ref{sec:eccentricity} we determine the limit average eccentricity of Fibonacci and Lucas cubes using related generating functions. In the subsequent section we then connect Fibonacci cubes with Fibonacci trees in a rather surprising way. Using a new labeling of the leaves of Fibonacci trees we prove that the eccentricity of a vertex of a given Fibonacci cube is equal to the depth of the associated leaf in the corresponding Fibonacci tree. Then, in Section~\ref{sec:weight}, we obtain the limit average fraction between the numbers of Fibonacci and Lucas strings with coordinates fixed to 0 and 1, respectively; see Theorem~\ref{thm:some-average} for the precise statement of the result. In the final section we first compute the limit average degree of these cubes. Then we introduce the hypercube density of a family of subgraphs of hypercubes and prove that it is equal to $(1-1/\sqrt 5)/\log_2\varphi$ for both Fibonacci cubes and Lucas cubes. We conclude the paper by demonstrating that the Cartesian product of graphs can be used to construct families of graphs with a fixed, non-zero hypercube density.   

The distance $d_G(u,v)$ between vertices $u$ and $v$ of a graph $G$ is the number of edges on a shortest $u,v$-path. The {\em eccentricity} ecc$_G(u)$ of $u\in V(G)$ is the maximum distance between $u$ and any other vertex of $G$. We will shortly write $d(u,v)$ and ecc$(u)$ when $G$ will be clear from the context. The {\em average eccentricity} and the {\em average degree} of a graph $G$ are respectively defined as:
\begin{eqnarray*}
\overline{{\rm ecc}}(G) & = & \frac{1}{|V(G)|}\sum_{u\in V(G)} {\rm ecc}(u)\,, \\ 
\overline{{\rm deg}}(G) & = & \frac{1}{|V(G)|}\sum_{u\in V(G)} {\rm deg}(u)\,.
\end{eqnarray*}

The vertex set of the \emph{$n$-cube} $Q_n$ is the set of all binary strings of length $n$, two vertices being adjacent if they differ in precisely one position.  
A {\em Fibonacci string} of length $n$ is a binary string $b_1\ldots b_n$ with $b_i\cdot b_{i+1}=0$ for $1\leq i<n$. The {\em Fibonacci cube} $\Gamma_n$ ($n\geq 1$) is the subgraph of $Q_n$ induced by the Fibonacci strings of length $n$. 
A Fibonacci string $b_1\ldots b_n$ is a {\em Lucas string} if in addition $b_1\cdot b_n = 0$ holds. The {\em Lucas cube} $\Lambda_n$ ($n\geq 1$) is the subgraph of $Q_n$ induced by the Lucas strings of length $n$. For convenience we also consider the empty string and set $\Gamma_0 = K_1 = \Lambda_0$. 

Let $\{F_n\}$ be the {\em Fibonacci numbers}: $F_0 = 0$, $F_1=1$, $F_{n} = F_{n-1} + F_{n-2}$ for $n\geq 2$. Recall that $\lim_{n\to \infty} F_{n+1}/F_n = \varphi$, where $\varphi = (1+\sqrt{5})/2$ is the {\em golden ratio}. More generally, if $k$ is a given integer, then $\lim_{n\to \infty} F_{n+k}/F_n = \varphi^k$.  
Let $\{L_n\}$ be {\em the Lucas numbers}: $L_0 = 2$, $L_1=1$, $L_{n} = L_{n-1} + L_{n-2}$ for $n\geq 2$. Recall finally that $|V(\Gamma_n)|= F_{n+2}$ for $n\ge 0$, and $|V(\Lambda_n)|= L_n = F_{n-1} + F_{n+1}$ for $n\ge 1$.

\section{Average eccentricity}
\label{sec:eccentricity}

In this section we determine the limit average eccentricity of Fibonacci and Lucas cubes. It is intuitively rather ``obvious" that in both cases the result should be the same, however, the proofs are somehow different. We begin with:

\begin{theorem}
\label{thm:ecc}
$$\lim_{n\to \infty} \frac{\overline{{\rm ecc}}(\Gamma_n)}{n} = \frac{5+\sqrt{5}}{10}\,.$$
\end{theorem}

\proof
Let $f_{n,k}$ be the number of vertices of $\Gamma_n$ with eccentricity $k$. It is proved in~\cite[Theorem 4.3]{camo-2012} that the corresponding generating function is
\begin{eqnarray}
\label{eq:gen-funct-fnk}
F(x,y) & = & \sum_{n,k\ge 0} f_{n,k} x^n y^k =  \frac{1+xy}{1-x(x+1)y}\,.
\end{eqnarray} 
If $e_n$ is the sum of the eccentricities of all vertices of $\Gamma_n$,
$$e_n = \sum_{x\in V(\Gamma_n)} {\rm ecc}(x)\,,$$
then
$$\frac{\partial F(x,y)}{\partial y}\bigg|_{y=1} = \sum_{n,k\ge 0} k f_{n,k} x^n = \sum_{n\ge 0} e_{n} x^n\,.$$
On the other hand, 
$$\frac{\partial F(x,y)}{\partial y}\bigg|_{y=1} = \frac{2x+x^2}{(1-x-x^2)^2}\,.$$
Since $\sum_{n\ge 0}F_nx^n = \frac{x}{1-x-x^2}$, we have 
\begin{eqnarray*}
\sum_{n\ge 0}F_{n+1}x^n & = & \frac{1}{1-x-x^2}, \\
\sum_{n\ge 0}nF_{n+1}x^n & = &  \frac{x+2x^2}{(1-x-x^2)^2},\ {\rm and}\\
\sum_{n\ge 0}nF_{n}x^n & = & \frac{x+x^3}{(1-x-x^2)^2}.
\end{eqnarray*}
Notice that 
$$\frac{2x+x^2}{(1-x-x^2)^2} = \frac{1}{5} \left( 3\frac{x}{1-x-x^2} + 4 \frac{x+2x^2}{(1-x-x^2)^2} + 3 \frac{x+x^3}{(1-x-x^2)^2}\right)\,.$$
Therefore, 
$$\sum_{n\ge 0}e_nx^n = \frac{1}{5} \left( 3\sum_{n\ge 0}F_{n}x^n + 4 \sum_{n\ge 0}nF_{n+1}x^n + 3 \sum_{n\ge 0}nF_{n}x^n\right)\,$$
and thus
$$e_n = \frac{3F_{n} + 4nF_{n+1} + 3nF_{n}}{5} = \frac{3F_{n} + nF_{n+1} + 3nF_{n+2}}{5}\,.$$
Therefore, 
$$\overline{{\rm ecc}}(\Gamma_n) = \frac{3F_{n} + nF_{n+1} + 3nF_{n+2}}{5F_{n+2}}\,.$$
We conclude that 
$$\lim_{n\to \infty} \frac{\overline{{\rm ecc}}(\Gamma_n)}{n} = \frac{3}{5} + \lim_{n\to \infty} \frac{1}{5} \frac{F_{n+1}}{F_{n+2}} = \frac{3}{5} + \frac{1}{5} \varphi^{-1} = \frac{5+\sqrt{5}}{10}\,.$$
\qed

Note that $(5+\sqrt{5})/10 \approx 0.7236$ which should be compared with the (trivial) fact that 
$$\lim_{n\to \infty} \frac{\overline{{\rm ecc}}(Q_n)}{n} = 1\,.$$

We next give the parallel result for Lucas cubes:
 
\begin{theorem}
\label{thm:eccluc}
$$\lim_{n\to \infty} \frac{\overline{{\rm ecc}}(\Lambda_n)}{n} =
\frac{5+\sqrt{5}}{10}\,.$$
\end{theorem}

\proof
The proof proceeds along the same lines as the proof of Theorem~\ref{thm:ecc}, but the computations are much different, hence we give a  sketch of the proof. Let $f'_{n,k}$ be the number of vertices of $\Lambda_n$ with eccentricity $k$. We start from the generating function of this sequence, obtained in~\cite[Theorem 5.16]{camo-2012}: 
\begin{eqnarray}
\label{eq:gen-funct-fpnk}
G(x,y) & = & \sum_{n,k\ge 0} f'_{n,k} x^n y^k =  \frac{1+x^2y}{1-xy-x^2y} +
\frac{1}{1+xy} - \frac{1-x}{1-x^2y}\,.
\end{eqnarray} 
Let $e'_n= \sum_{x\in V(\Lambda_n)} {\rm ecc}(x)$ be the sum of the eccentricities of all vertices of $\Lambda_n$. We deduce from~\eqref{eq:gen-funct-fpnk} that the generating function of the sequence $\{e'_n\}$ is
$$\sum_{n\ge 0} e'_{n} x^n = \frac{\partial G(x,y)}{\partial y}\bigg|_{y=1}
=\frac{x+2x^2}{(1-x-x^2)^2}\ - \frac{x}{(1+x)^2} - \frac{x^2}{(1+x)(1-x^2)} \,.$$
The first term is the generating function of $n\;  F_{n+1}$, and developing the other terms we obtain
 $$e'_n = n\;  F_{n+1} \; +(-1)^n n\;  +(-1)^{n+1}\left\lfloor\frac{n}{2}
\right\rfloor.$$
Since $|V(\Lambda_n)|=F_{n-1}+F_{n+1}$, we conclude that 
$$\lim_{n\to \infty} \frac{\overline{{\rm ecc}}(\Lambda_n)}{n} = \lim_{n\to \infty}
\frac{F_{n+1}}{F_{n-1}+F_{n+1}}  = \frac{\varphi^{2}}{1+\varphi^{2}} =
\frac{5+\sqrt{5}}{10}\,.$$
\qed

\section{Fibonacci trees and eccentricity}
\label{sec:eccentricitysequence}

In the previous section we considered the sequence $e_n$, where $e_n = \sum_{x\in V(\Gamma_n)} {\rm ecc}(x)$. The sequence $e_{n+1}$ starts with $2,5,12,25,50,96,\ldots$ This is also the start of the sequence~\cite[Sequence A067331]{sloane} described as the sum of the depth of leaves in Fibonacci trees. The two sequences indeed coincide since they have the same generating function but the connection seems mysterious. In this section we give a bijective proof that the two sequences coincide and along the way propose a new labeling of the Fibonacci trees.  

Fibonacci trees have been introduced in computer science in the context of efficient search algorithms~\cite{ho-1982,knuth-1998,wa-2007}. They are binary trees defined recursively as follows:
\begin{itemize}
\item $T_0$ and $T_1$ are trees with a single vertex---the root.
\item $T_n$, $n\geq 2$, is the rooted tree whose left subtree is  $T_{n-1}$ and whose right subtree is  $T_{n-2}$.
\end{itemize}
Clearly, for any $n$ the number of leaves of $T_n$ is $F_{n+1}$. 

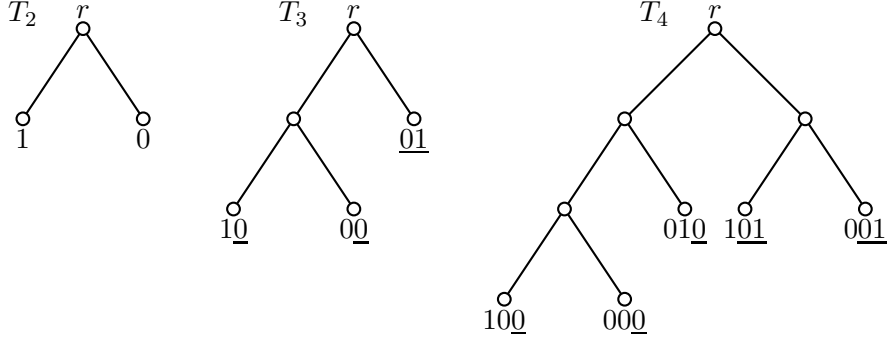
\begin{figure}[ht!]
\begin{center}
\begin{tikzpicture}[scale=0.4,style=thick]
\def\vr{6pt} 
\path (2,9) coordinate (r1); 
\path (0,6) coordinate (a1); 
\path (4,6) coordinate (a2); 
\path (11,9) coordinate (r2); 
\path (9,6) coordinate (b1); 
\path (13,6) coordinate (b2); 
\path (7,3) coordinate (b3); 
\path (11,3) coordinate (b4); 
\path (23,9) coordinate (r3); 
\path (20,6) coordinate (c1); 
\path (26,6) coordinate (c2); 
\path (18,3) coordinate (c3); 
\path (22,3) coordinate (c4); 
\path (24,3) coordinate (c5); 
\path (28,3) coordinate (c6); 
\path (16,0) coordinate (c7); 
\path (20,0) coordinate (c8); 
\draw (a1) -- (r1) -- (a2);
\draw (b3) -- (b1) -- (r2) --(b2);
\draw (b1) -- (b4);
\draw (c7) -- (c3) -- (c1) --(r3) -- (c2) -- (c6);
\draw (c3) -- (c8); \draw (c1) -- (c4); \draw (c2) -- (c5);
\draw (r1)  [fill=white] circle (\vr);
\draw (a1)  [fill=white] circle (\vr);
\draw (a2)  [fill=white] circle (\vr);
\draw (r2)  [fill=white] circle (\vr);
\draw (b1)  [fill=white] circle (\vr);
\draw (b2)  [fill=white] circle (\vr);
\draw (b3)  [fill=white] circle (\vr);
\draw (b4)  [fill=white] circle (\vr);
\draw (r3)  [fill=white] circle (\vr);
\draw (c1)  [fill=white] circle (\vr);
\draw (c2)  [fill=white] circle (\vr);
\draw (c3)  [fill=white] circle (\vr);
\draw (c4)  [fill=white] circle (\vr);
\draw (c5)  [fill=white] circle (\vr);
\draw (c6)  [fill=white] circle (\vr);
\draw (c7)  [fill=white] circle (\vr);
\draw (c8)  [fill=white] circle (\vr);

\draw [above] (r1) node {$r$};
\draw [above] (r2) node {$r$};
\draw [above] (r3) node {$r$};
\draw [below] (a1) node {$1$};
\draw [below] (a2) node {$0$};
\draw [below] (b3) node {$1\underline{0}$};
\draw [below] (b4) node {$0\underline{0}$};
\draw [below] (b2) node {$\underline{01}$};
\draw [below] (c7) node {$10\underline{0}$};
\draw [below] (c8) node {$00\underline{0}$};
\draw [below] (c4) node {$01\underline{0}$};
\draw [below] (c5) node {$1\underline{01}$};
\draw [below] (c6) node {$0\underline{01}$};
\draw (0,9.5) node {$T_2$};
\draw (9,9.5) node {$T_3$};
\draw (21,9.5) node {$T_4$};
\end{tikzpicture}
\end{center}
\caption{Fibonacci trees equipped with the standard labeling}
\label{fig:standard}
\end{figure}

We can recursively construct a labeling of the leaves of $T_n$ with Fibonacci strings of length $n-1$ as follows. Let ${\cal F}_n$ be the set of Fibonacci strings of length $n$. Let ${\cal F}_n^0$ and ${\cal F}_n^1$ be the sets of Fibonacci strings ending with $0$ and $1$, respectively. We then have, for $n\geq 2$, ${\cal F}_n = {\cal F}_n^0 \uplus {\cal F}_n^1 = \{s0;s\in {\cal F}_{n-1}\}  \uplus \{s01;s\in {\cal F}_{n-2}\}$, where $\uplus$ is the disjoint union of sets. First label $T_1$ and $T_2$ and assume $n\geq 3$.  We append $01$ to the labels of the right leaves already labeled as leaves of  $T_{n-2}$, and  $0$ to the labels of the left leaves already labeled as leaves of  $T_{n-1}$. In Fig.~\ref{fig:standard} the construction is presented for the first three non-trivial Fibonacci trees, where the currently attached strings are underlined.

The described standard labeling of Fibonacci trees does not respect the equality between depth and eccentricity. For example, the depth of the leaf of $T_3$ labeled $01$ is $1$, but ${\rm ecc}_{\Gamma_2}(01)=2$. Nevertheless, the sum of the depths of leaves of $T_3$ is $5$ like the sum of the eccentricities of vertices of $\Gamma_2$. We next construct a new labeling that respects the equality vertex by vertex. 

\begin{figure}[ht!]
\begin{center}
\begin{tikzpicture}[scale=0.4,style=thick]
\def\vr{6pt} 
\path (2,9) coordinate (r1); 
\path (0,6) coordinate (a1); 
\path (4,6) coordinate (a2); 
\path (11,9) coordinate (r2); 
\path (9,6) coordinate (b1); 
\path (13,6) coordinate (b2); 
\path (7,3) coordinate (b3); 
\path (11,3) coordinate (b4); 
\path (23,9) coordinate (r3); 
\path (20,6) coordinate (c1); 
\path (26,6) coordinate (c2); 
\path (18,3) coordinate (c3); 
\path (22,3) coordinate (c4); 
\path (24,3) coordinate (c5); 
\path (28,3) coordinate (c6); 
\path (16,0) coordinate (c7); 
\path (20,0) coordinate (c8); 
\draw (a1) -- (r1) -- (a2);
\draw (b3) -- (b1) -- (r2) --(b2);
\draw (b1) -- (b4);
\draw (c7) -- (c3) -- (c1) --(r3) -- (c2) -- (c6);
\draw (c3) -- (c8); \draw (c1) -- (c4); \draw (c2) -- (c5);
\draw (r1)  [fill=white] circle (\vr);
\draw (a1)  [fill=white] circle (\vr);
\draw (a2)  [fill=white] circle (\vr);
\draw (r2)  [fill=white] circle (\vr);
\draw (b1)  [fill=white] circle (\vr);
\draw (b2)  [fill=white] circle (\vr);
\draw (b3)  [fill=white] circle (\vr);
\draw (b4)  [fill=white] circle (\vr);
\draw (r3)  [fill=white] circle (\vr);
\draw (c1)  [fill=white] circle (\vr);
\draw (c2)  [fill=white] circle (\vr);
\draw (c3)  [fill=white] circle (\vr);
\draw (c4)  [fill=white] circle (\vr);
\draw (c5)  [fill=white] circle (\vr);
\draw (c6)  [fill=white] circle (\vr);
\draw (c7)  [fill=white] circle (\vr);
\draw (c8)  [fill=white] circle (\vr);

\draw [above] (r1) node {$r$};
\draw [above] (r2) node {$r$};
\draw [above] (r3) node {$r$};
\draw [below] (a1) node {$1$};
\draw [below] (a2) node {$0$};
\draw [below] (b3) node {$1\underline{0}$};
\draw [below] (b4) node {$0\underline{1}$};
\draw [below] (b2) node {$\underline{00}$};
\draw [below] (c7) node {$10\underline{1}$};
\draw [below] (c8) node {$01\underline{0}$};
\draw [below] (c4) node {$00\underline{1}$};
\draw [below] (c5) node {$1\underline{00}$};
\draw [below] (c6) node {$0\underline{00}$};
\draw (0,9.5) node {$T_2$};
\draw (9,9.5) node {$T_3$};
\draw (21,9.5) node {$T_4$};
\end{tikzpicture}
\end{center}
\caption{Fibonacci trees equipped with the labeling $\theta$}
\label{fig:new}
\end{figure}
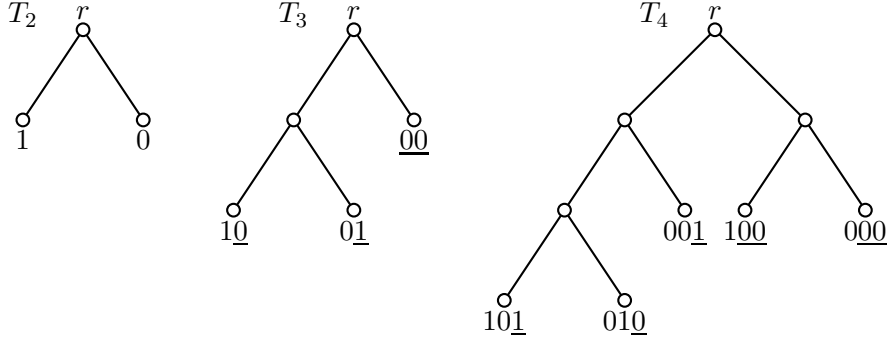

Notice first that, for $n\geq 2$, ${\cal F}_{n} = \{s00;s\in {\cal F}_{n-2}\}\uplus \{s0;s\in {\cal F}_{n-1}^1\}\uplus \{s1;s\in {\cal F}_{n-1}^0\}$. Label $T_1$ with the empty string and $T_2$ according to Fig.~\ref{fig:new}. Assume that $n\geq 3$ and that the leaves of $T_{n-1}$ and $T_{n-2}$ were already labeled. We  append $00$ to the label of the right leaves. For the left leaves append $0$ to the label of a leaf with a label ending with $1$; otherwise append label $1$. Let $\theta$ denote this  labeling of the leaves of $T_n$ by vertices of $\Gamma_{n-1}$. It is shown in Fig.~\ref{fig:new}, again for the first three non-trivial Fibonacci trees.

The main result of this section now reads as follows:

\begin{theorem}
\label{thm:Fibo-tree}
Let $n\ge 1$ and let $u\in V(\Gamma_n)$. Then 
$${\rm ecc}_{\Gamma_{n}}(u) = {\rm depth}_{T_{n+1}}(\theta^{-1}(u))\,.$$
\end{theorem}

\proof
We proceed by induction on $n$, the cases $n=1,2$ being trivial. Let $n\ge 3$ and let $u\in V(\Gamma_n)$. Consider the following three cases. 

Suppose first that $u=v00$, where $v\in {\cal F}_{n-2}$. Then we claim that ${\rm ecc}_{\Gamma_{n}}(u) = {\rm ecc}_{\Gamma_{n-2}}(v) + 1$. Let $u'$ be a vertex of $\Gamma_n$ with $d_{\Gamma_n}(u,u') = {\rm ecc}(u)$. Let $u' = wab$, where $w\in {\cal F}_{n-2}$. Since $ab\neq 11$, $d_{\Gamma_n}(u,u') \neq d_{\Gamma_{n-2}}(v,w) +2$.  Then 
$${\rm ecc}_{\Gamma_{n}}(u) = d_{\Gamma_n}(u,u') \le d_{\Gamma_{n-2}}(v,w) + 1 \le {\rm ecc}(v)_{\Gamma_{n-2}} + 1\,.$$  
Conversely, let $w\in {\cal F}_{n-2}$, such that $d_{\Gamma_{n-2}}(v,w) = {\rm ecc}(v)$. Then $w01\in {\cal F}_{n}$ and $d_{\Gamma_{n}}(u,w01) = {\rm ecc}_{\Gamma_{n-2}}(v)+1$. If follows that ${\rm ecc}_{\Gamma_{n}}(u) \ge {\rm ecc}(v)_{\Gamma_{n-2}} + 1$. This proves the claim. 

Suppose next that $u=v1$, where $v\in {\cal F}_{n-1}^{0}$. Now we claim that ${\rm ecc}_{\Gamma_{n}}(u) = {\rm ecc}_{\Gamma_{n-1}}(v) + 1$. The inequality ${\rm ecc}_{\Gamma_{n}}(u) \le {\rm ecc}_{\Gamma_{n-1}}(v) + 1$ follows by an argument similar as in the first case.  
Conversely, let $w\in {\cal F}_{n-1}$, such that $d_{\Gamma_{n-1}}(v,w) = {\rm ecc}(v)$. Then $w0\in {\cal F}_{n}$ and $d_{\Gamma_{n}}(u,w0) = {\rm ecc}_{\Gamma_{n-1}}(v)+1$. If follows that ${\rm ecc}_{\Gamma_{n}}(u) \ge {\rm ecc}_{\Gamma_{n-1}}(v) + 1$. 

Suppose finally that $u=v0$, where $v\in {\cal F}_{n-1}^{1}$. We claim again that ${\rm ecc}_{\Gamma_{n}}(u) = {\rm ecc}_{\Gamma_{n-1}}(v) + 1$. Again, the inequality ${\rm ecc}_{\Gamma_{n}}(u) \le {\rm ecc}_{\Gamma_{n-1}}(v) + 1$ follows as above.  
Conversely, let $w\in {\cal F}_{n-1}$, such that $d_{\Gamma_{n-1}}(v,w) = {\rm ecc}_{\Gamma_{n-1}}(v)$. Then $w$ ends with 0, because otherwise $w$ would not be an eccentric vertex of $v$. Indeed, if $w$ ended with 1, then the word $w'$ obtained from $w$ by changing its last bit to 0 would lie in ${\cal F}_{n-1}$ and hence $d_{\Gamma_{n-1}}(w',v) > d_{\Gamma_{n-1}}(w,v) = {\rm ecc}_{\Gamma_{n-1}}(v)$, a contradiction. 
It follows that $w1\in {\cal F}_{n}$ and $d_{\Gamma_{n}}(u,w1) = d_{\Gamma_{n-1}}(w,v) + 1 = {\rm ecc}_{\Gamma_{n-1}}(v)+1$. 

We have thus proved that for any $u\in {\cal F}_{n}$, ${\rm ecc}_{\Gamma_{n}}(u)$ increases by 1 with respect to the word $v$ to which a suffix has been added to obtain $u$. By the construction of $T_n$ and by the induction hypothesis, 
$${\rm depth}_{T_{n+1}}(\theta^{-1}(u)) = {\rm depth}_{T_{n-1}}(\theta^{-1}(v))+1 = {\rm ecc}_{\Gamma_{n-2}}(v)+1={\rm ecc}_{\Gamma_{n}}(u)$$ 
in the first case, and 
$${\rm depth}_{T_{n+1}}(\theta^{-1}(u)) = {\rm depth}_{T_{n}}(\theta^{-1}(v))+1 = {\rm ecc}_{\Gamma_{n-1}}(v)+1={\rm ecc}_{\Gamma_{n}}(u)$$ 
in the last two cases. 
\qed

\section{Average fractional weights}
\label{sec:weight}

Let $G$ be a subgraph of $Q_n$, so that the vertices of $G$ are binary strings of length $n$. Then for $i=1,\ldots, n$ and $\chi =0,1$, let  
$$W_{(i,\chi)}(G)=\{ u=u_1\ldots u_n\in V(G)\ |\ u_i=\chi\}\,.$$
In this section we prove the following result which might be of independent interest: 

\begin{theorem}
\label{thm:some-average}
$$\lim_{n\to \infty} \frac{1}{n}  \sum_{i=1}^{n}
  \frac{|W_{(i,0)}(\Gamma_n)|} {|W_{(i,1)}(\Gamma_n)|} = 
  \lim_{n\to \infty} \frac{1}{n}  \sum_{i=1}^{n}
  \frac{|W_{(i,0)}(\Lambda_n)|} {|W_{(i,1)}(\Lambda_n)|} 
= \varphi^2\,.$$
Moreover, for every $i$ between 1 and $n$, we have
$$\lim_{n\to \infty}    \frac{|W_{(i,0)}(\Lambda_n)|}{|W_{(i,1)}(\Lambda_n)|} 
= \varphi^2\,.$$
\end{theorem}

\noindent
To prove it, we will make use of the following result, see~\cite[Exercise 3.9.13]{monier-1990}: 

\begin{lemma}
\label{lem:classic}
If $\{a_n\}$ is a convergent complex sequence with limit $\ell$, then 
$$\lim_{n\to \infty} \frac{1}{n}\sum_{i=1}^n a_ia_{n+1-i} = \ell^2\,.$$
\end{lemma}

\proof (of Theorem~\ref{thm:some-average})
Set $w_{i,0}^{(n)} = |W_{(i,0)}(\Gamma_n)|$, $w_{i,1}^{(n)} = |W_{(i,1)}(\Gamma_n)|$, and $x_i^{(n)} = w_{i,0}^{(n)} / w_{i,1}^{(n)}$. Notice first that the vertices of $W_{(i,0)}(\Gamma_n)$ are the strings $u0v$ where $u$ and $v$ are arbitrary vertices of $V(\Gamma_{i-1})$ and $V(\Gamma_{n-i})$, respectively. Similarly the vertices of $W_{(i,1)}(\Gamma_n)$ are the strings $u010v$, where $u$ and $v$ are arbitrary vertices of $V(\Gamma_{i-2})$ and $V(\Gamma_{n-i-1})$, respectively.
Then, having in mind that $|V(\Gamma_n)|=F_{n+2}$, we have 
$$x_i^{(n)} = \frac{F_{i+1}\cdot F_{n-i+2}}{F_{i}\cdot F_{n-i+1}}\,.$$
Setting $a_n = F_{n+1}/F_n$, recalling that $\{a_n\} \to \varphi$, and using Lemma~\ref{lem:classic}, we get
$$\lim_{n\to \infty} \frac{1}{n}\sum_{i=1}^n x_i^{(n)} = 
\lim_{n\to \infty} \frac{1}{n}\sum_{i=1}^n a_i a_{n-i+1} = \varphi^2\,,$$ 
which proves the result for Fibonacci cubes. 

For Lucas cubes we have $|W_{(i,0)}(\Lambda_n)| = F_{n+1}$ and $|W_{(i,1)}(\Lambda_n)| = F_{n-1}$ for all $1\le i\le n$. This can be seen by considering Lucas strings not as linear orderings, but rather as cyclic orderings of $0$'s and $1$'s with no consecutive $1$'s. Then removing a $0$ in the $i$-th position of the cycle gives a bijection between
$W_{(i,0)}(\Lambda_n)$ and $V(\Gamma_{n-1})$, while removing a segment $010$ centered at the $i$-th position of the cycle gives a bijection between $W_{(i,1)}(\Lambda_n)$ and $V(\Gamma_{n-3})$. Consequently, if $w_{i,0}^{(n)} = |W_{(i,0)}(\Lambda_n)|$, $w_{i,1}^{(n)} = |W_{(i,1)}(\Lambda_n)|$, and $x_i^{(n)} = w_{i,0}^{(n)} / w_{i,1}^{(n)}$, then $x_i^{(n)} = F_{n+1}/F_{n-1}$. It follows that $\lim_{n\to \infty} x_i^{(n)} = \varphi^2$. Then the classical Ces\`{a}ro Means Theorem (it asserts that if $\lim_{n\to \infty} a_n = \ell$, then $\lim_{n\to \infty} \frac{1}{n}\sum_{i=1}^n a_i = \ell$ as well) can be applied instead of Lemma~\ref{lem:classic}. 
\qed

\section{Average degree and density}
\label{sec:degree}

In this section  we first compute the limit average degree of the considered graphs. The result then motivates us to introduce the hypercube density, to determine it for the cubes, and to show that the Cartesian product of graphs is useful in this context. 

\begin{theorem}
\label{thm:degree}
$$\lim_{n\to \infty} \frac{\overline{{\rm deg}}(\Gamma_n)}{n} = \lim_{n\to \infty} \frac{\overline{{\rm deg}}(\Lambda_n)}{n} = \frac{5-\sqrt{5}}{5}\,.$$
\end{theorem}

\proof
It was proved in~\cite{mupe-2001} that $|E(\Gamma_n)| = (nF_{n+1} + 2(n+1)F_n)/5$, hence 
$$\overline{{\rm deg}}(\Gamma_n) = \frac{1}{|V(\Gamma_n)|}2\,|E(\Gamma_n)| = 
\frac{2}{5\,F_{n+2}}(nF_{n+1} + 2(n+1)F_n)\,.$$
Therefore, 
$$\lim_{n\to \infty} \frac{\overline{{\rm deg}}(\Gamma_n)}{n} = \lim_{n\to \infty} \frac{2}{5}\left(\frac{F_{n+1}}{F_{n+2}} + \frac{2F_{n}}{F_{n+2}} + \frac{2}{n}\frac{F_{n}}{F_{n+2}}\right) = \frac{2}{5}\left(\varphi^{-1} + 2\varphi^{-2}\right) = \frac{5-\sqrt{5}}{5}\,.$$ 

For the Lucas cubes we recall from~\cite[p.1322]{klmope-2011} that $|E(\Lambda_n)| = nF_{n-1}$. Hence $\overline{{\rm deg}}(\Lambda_n) = 2nF_{n-1}/(F_{n-1}+F_{n+1})$ and 
$$\lim_{n\to \infty} \frac{\overline{{\rm deg}}(\Lambda_n)}{n} 
= \lim_{n\to \infty} \frac{2F_{n-1}}{F_{n-1} + F_{n+1}} 
= \frac{2}{1+\varphi^{2}} = \frac{5-\sqrt{5}}{5}\,.$$
\qed

Graham~\cite{gr-1970} proved the following fundamental property of subgraphs of hypercubes (see~\cite[Lemma 3.2]{haim-2011} for an alternative proof of it): 

\begin{lemma} [Density Lemma]
\label{lem:Graham}
\label{cha12-15}
Let $G$ be a subgraph of a hypercube. Then
$$|E(G)| \leq {1\over 2} |V(G)|\cdot \log_2 |V(G)|\,.$$
Moreover, equality holds if and only if $G$ is a hypercube.
\end{lemma}

The lemma has important consequences, in particular for fast recognition 
algorithms for classes of subgraphs of hypercubes; see~\cite{tave-2007} for the case of Fibonacci cubes. The Density Lemma also asserts that hypercubes have the largest density among all subgraphs of hypercubes. We therefore introduce the following concept. 

If $G$ is a subgraph of a hypercube, then let 
$$\rho(G) = \frac{\overline{{\rm deg}}(G)}{\log_2 |V(G)|}\,.$$
Let ${\cal G} = \{G_k\}_{k\ge 1}$ be an {\em increasing family} of subgraphs of hypercubes, that is, a family with $|V(G_{n+1})| > |V(G_{n})|$ for $n\ge 1$. Then the {\em hypercube density} of ${\cal G}$ is 
$$\rho({\cal G}) = \limsup_{k\to \infty} \rho(G_k)\,.$$
By the Density Lemma, $0\le \rho(\{G_k\}) \le 1$ holds for any family $\{G_k\}$ and 
$$\rho(\{Q_k\}) = 1\,.$$ 
For Fibonacci cubes and Lucas cubes we have:

\begin{corollary}
\label{coro:densityFibonacciLucas}
$$\rho(\{\Gamma_n\}) = \rho(\{\Lambda_n\}) = \frac{5-\sqrt5}{5\log_2 \varphi}\,.$$
\end{corollary}

\proof
The result easily follows from Theorem~\ref{thm:degree} together with the facts that $F_n \sim \varphi^n/\sqrt 5$ and that $\Lambda_n \sim \varphi^n$. 
\qed

Hence $\rho(\{\Gamma_n\})\approx 0.7962$ which is in particular interesting because the hypercube density of many other important families of hypercube subgraphs is 0. For instance, if $WB_k$ denotes the bipartite wheel with $k$ spokes, then an easy calculation shows that $\rho(\{WB_k\}) = 0$. (Cf.~\cite{bach-1996} for the role of bipartite wheels among subgraphs of hypercubes.) For another example consider the subdivision $S(K_k)$ of the complete graph $K_k$, that is, the graph obtained from $K_k$ by subdividing each of its edges precisely once. These graphs embed isometrically into hypercubes (cf.~\cite{klli-2003}) and 
$$\rho(\{S(K_k)\}) = 
\lim_{k\to \infty} \rho(S(K_k)) = 
\lim_{k\to \infty} \frac{2 \left(2 {k\choose 2}\right)} {\left(k + {k\choose 2}\right) \log_2\left(k + {k\choose 2}\right)} = 
0\,.$$

Families of graphs with bounded degree also have density equal to 0. More precisely: 

\begin{proposition}
\label{prp:bounded-degree}
Let $\{{G_k}\}$ be an increasing family of hypercube subgraphs. If there exists a constant $M$ such that $\Delta(G_k) \le M$ for any $k\ge 1$, then $\rho(\{G_k\}) = 0$. 
\end{proposition}

\proof
For $k\ge 1$ set $n_k=|V(G_k)|$ and $m_k=|E(G_k)|$. Then from the Handshaking Lemma it follows 
that $2m_k = \sum_{u\in V(G_k)} {\rm deg}(u) \le n_k\, M$ and hence $2m_k/n_k \le M$. The assertion 
then follows because $\log_2n_k \to \infty$. 
\qed

On the other hand, the Cartesian product of graphs can be used to obtain families with positive hypercube density. If $G^k$ denotes the $k$-tuple Cartesian product of $G$, then we have: 

\begin{proposition}
\label{prp:Cart-powers}
Let $G$ be a hypercube subgraph with $\rho(G) = c$. Then $\rho( \{G^k\} ) = c$. 
\end{proposition}
 
\proof
Let $n=|V(G)|$ and $m=|E(G)|$. Then it is easily shown by induction that 
$|V(G^k)| = n^k$ and $|E(G^k)| = kn^{k-1}m$. Consequently, 
$$\rho(\{ G^k \} ) = \lim_{k\to \infty} \rho(G^k) = \lim_{k\to \infty} \frac{2kn^{k-1}m}{n^k \log_2n^k}  
= \frac{2m}{n\log_2n} = \rho(G)\,.$$
\qed

Suppose that a family $\{G_k\}$ is given with $\rho(\{{G_k}\}) = c$. Then we can again use the Cartesian product of graphs to obtain an infinite number of families with the same density. 

\begin{proposition}
\label{prp:products}
Let $\rho(\{{G_k}\}) = c$ and let $G$ be a fixed subgraph of a hypercube. Then $\rho( \{G_k\cp G \} ) = c$. 
\end{proposition}

\proof
Let $n=|V(G)|$, $m=|E(G)|$, and for $k\ge 1$, let $n_k=|V(G_k)|$ and $m_k=|E(G_k)|$. Then we have
\begin{eqnarray*}
\rho(\{ G_k\cp G \} ) & = & \lim_{k\to \infty} \rho(G_k\cp G) \\
& = & \lim_{k\to \infty} \frac{2(n\,m_k + n_k\, m)}{n_k\,n \log_2(n_k\,n)} \\
& = & \lim_{k\to \infty} \frac{2m_k}{n_k\log_2(n_k\,n)} +
      \lim_{k\to \infty} \frac{2m}{n\log_2(n_k\,n)} \\
& = & \lim_{k\to \infty} \left( 
\frac{2m_k}{n_k\log_2n_k}\cdot \frac{\log_2n_k}{\log_2n_k + \log_2n}
\right) \\
& = & \lim_{k\to \infty} \frac{2m_k}{n_k\log_2n_k} = \rho(\{{G_k}\})\,.
\end{eqnarray*}
\qed

We conclude the paper with the following question proposed to us by one the referees: 

\begin{problem}
Can every real value in $[0,1]$ be obtained as hypercube density of an increasing family? 
\end{problem}

\section*{Acknowledgments}

This work was supported in part by the Proteus project 
BI-FR/12-139-PROTEUS-008, by ARRS Slovenia under the grant
P1-0297, and within the EUROCORES Programme EUROGIGA/GReGAS of the
European Science Foundation. 

We thank the referees for careful reading of the manuscript and useful remarks and suggestions.


\begin{thebibliography}{99}

\bibitem{bach-1996}
  H.-J.~Bandelt, V.~Chepoi,
  Graphs of acyclic cubical complexes,
  European J. Combin. 17 (1996) 113--120. 

\bibitem{camo-2012}
  A.~Castro, M.~Mollard,
  The eccentricity sequences of Fibonacci and Lucas cubes,
  Discrete Math. 312 (2012) 1025--1037.

\bibitem{dago-2004}
  P.~Dankelmann, W.~Goddard, C.S.~Swart, 
  The average eccentricity of a graph and its subgraphs, 
  Util. Math. 65 (2004) 41--51.
  
\bibitem{gr-1970}
  R.~L.~Graham, 
  On primitive graphs and optimal vertex assignments, 
  Ann. New York Acad. Sci. 175 (1970) 170--186.  

\bibitem{haim-2011}
  R.~Hammack, W.~Imrich, S.~Klav\v zar, 
  Handbook of Product Graphs, Second Edition,
  CRC Press, Boca Raton, FL, 2011. 

\bibitem{hipa-2012}
  A.M.~Hinz, D.~Parisse,
  The average eccentricity of Sierpi\'nski graphs,
  Graphs Combin. 28 (2012) 671--686.
   
\bibitem{ho-1982}
  Y.~Horibe,
  An entropy view of Fibonacci trees,
  Fibonacci Quart. 20 (1982) 168--178. 
  
\bibitem{hsu-1993}
  W.-J. Hsu,
  Fibonacci cubes{---}a new interconnection technology,
  IEEE Trans. Parallel Distrib. Syst. 4 (1993) 3--12.

\bibitem{il-2012}
  A.~Ili\'c,
  On the extremal properties of the average eccentricity,
  Comp. Math. Appl. 64 (2012) 2877--2885.

\bibitem{kl-2013}
  S.~Klav\v{z}ar, 
  Structure of Fibonacci cubes: a survey, 
  J. Comb. Optim. 25 (2013) 505--522.  

\bibitem{klli-2003}  
  S.~Klav\v zar, A.~Lipovec, 
  Partial cubes as subdivision graphs and as generalized Petersen graphs,
  Discrete Math. 263 (2003) 157--165. 

\bibitem{klmo-2012}
  S.~Klav\v{z}ar, M.~Mollard,
  Wiener index and Hosoya polynomial of Fibonacci and Lucas cubes,
  MATCH Commun. Math. Comput. Chem. 68 (2012) 311--324. 
  
\bibitem{klmope-2011}
  S.~Klav\v{z}ar, M.~Mollard, M.~Petkov\v sek,
  The degree sequence of Fibonacci and Lucas cubes,
  Discrete Math. 311 (2011) 1310--1322.

\bibitem{knuth-1998} 
  D.E.~Knuth,
  The Art of Computer Programming, Vol. 3, 2nd edition, 
  Addison-Wesley, Reading, MA, 1998, p. 417.
  
\bibitem{monier-1990}
  J.-M.~Monier,
  Analyse, Tome 1,
  Dunod, 1990.    

\bibitem{mupe-2001}
  E.~Munarini, C.~Perelli Cippo, N.~Zagaglia Salvi,
  On the Lucas cubes,
  Fibonacci Quart. 39 (2001) 12--21.

\bibitem{ra-2013}
  M.~Ramras,
  Routing linear permutations on Fibonacci and Lucas cubes,
  manuscript, arXiv:1207.1518v1 [math.CO].

\bibitem{sloane}
  N.J.A.~Sloane,
  The {O}n-{L}ine {E}ncyclopedia of {I}nteger {S}equences,
  Published electronically at http://oeis.org, 2013.

\bibitem{tave-2007}
  A.~Taranenko, A.~Vesel, 
  Fast recognition of Fibonacci cubes,
  Algorithmica 49 (2007) 81--93.

\bibitem{wa-2007}
  S.G.~Wagner,
  The Fibonacci number of Fibonacci trees and a related family of polynomial recurrence systems,
  Fibonacci Quart. 45 (2007) 247--253. 

\bibitem{zibe-2012}
  P.~\v{Z}igert Pleter\v{s}ek, M.~Berli\v{c},
  Resonance graphs of armchair nanotubes cyclic polypyrenes and amalgams of Lucas cubes, 
  MATCH Commun. Math. Comput. Chem. 70 (2013) 533--543.

\bibitem{zhou-2009}
  H.~Zhang, L.~Ou, H.~Yao,
  Fibonacci-like cubes as {$Z$}-transformation graphs,
  Discrete Math. 309 (2009) 1284--1293. 

\end{thebibliography}
\end{document}